\documentclass[]{amsart}
\usepackage{amsmath, amsthm, amscd, amsfonts, amssymb, graphicx, color}
\usepackage[bookmarksnumbered, colorlinks, plainpages]{hyperref}

\theoremstyle{definition}

\theoremstyle{remark}

\numberwithin{equation}{section}

\begin{document}
\setcounter{page}{1}
\begin{center}
{\bf THE TOPOLOGICAL CENTERS OF MODULE ACTIONS}
\end{center}

\title[]{}

\author[]{KAZEM HAGHNEJAD AZAR  AND ABDOLHAMID RIAZI}

\address{}

\dedicatory{}



\begin{abstract}
In this paper, we will study the topological centers of $n-th$ dual of Banach $A-module$ and we extend some propositions from Lau and $\ddot{U}lger$ into $n-th$ dual of Banach $A-modules$. Let $B$ be a Banach  $A-bimodule$. By using some new conditions, we show that ${ {Z}^\ell}_{A^{(n)}}(B^{(n)})=B^{(n)}$ and ${ {Z}^\ell}_{B^{(n)}}(A^{(n)})=A^{(n)}$. We also have some conclusions in  dual groups.
\end{abstract} \maketitle

\section{\bf  Preliminaries and
Introduction }

\noindent Throughout  this paper, $A$ is  a Banach algebra and $A^*$,
$A^{**}$, respectively, are the first and second dual of $A$. Recall that  a left approximate identity $(=LAI)$ [resp. right
approximate identity $(=RAI)$]
in Banach algebra $A$ is a net $(e_{\alpha})_{{\alpha}\in I}$ in $A$ such that   $e_{\alpha}a\longrightarrow a$ [resp. $ae_{\alpha}\longrightarrow a$]. We
say that a  net $(e_{\alpha})_{{\alpha}\in I}\subseteq A$ is a
approximate identity $(=AI)$ for $A$ if it is $LAI$ and $RAI$ for $A$. If $(e_{\alpha})_{{\alpha}\in I}$ in $A$ is bounded and $AI$ for $A$, then we say that $(e_{\alpha})_{{\alpha}\in I}$ is a bounded approximate identity  ($=BAI$) for $A$.
 For $a\in A$
 and $a^\prime\in A^*$, we denote by $a^\prime a$
 and $a a^\prime$ respectively, the functionals on $A^*$ defined by $<a^\prime a,b>=<a^\prime,ab>=a^\prime(ab)$ and $<a a^\prime,b>=<a^\prime,ba>=a^\prime(ba)$ for all $b\in A$.
 The Banach algebra $A$ is embedded in its second dual via the identification
 $<a,a^\prime>$ - $<a^\prime,a>$ for every $a\in
A$ and $a^\prime\in
A^*$.
  We denote the set   $\{a^\prime a:~a\in A~ and ~a^\prime\in
  A^*\}$ and
  $\{a a^\prime:~a\in A ~and ~a^\prime\in A^*\}$ by $A^*A$ and $AA^*$, respectively, clearly these two sets are subsets of $A^*$.\\
   Let $A$ have a $BAI$. If the
equality $A^*A=A^*,~~(AA^*=A^*)$ holds, then we say that $A^*$
factors on the left (right). If both equalities $A^*A=AA^*=A^*$
hold, then we say
that $A^*$  factors on both sides.\\
 The extension of bilinear maps on normed space and the concept of regularity of bilinear maps were studied by [1, 2, 3, 6, 8, 14]. We start by recalling these definitions as follows.\\
 Let $X,Y,Z$ be normed spaces and $m:X\times Y\rightarrow Z$ be a bounded bilinear mapping. Arens in [1] offers two natural extensions $m^{***}$ and $m^{t***t}$ of $m$ from $X^{**}\times Y^{**}$ into $Z^{**}$ as following\\
1. $m^*:Z^*\times X\rightarrow Y^*$,~~~~~given by~~~$<m^*(z^\prime,x),y>=<z^\prime, m(x,y)>$ ~where $x\in X$, $y\in Y$, $z^\prime\in Z^*$,\\
 2. $m^{**}:Y^{**}\times Z^{*}\rightarrow X^*$,~~given by $<m^{**}(y^{\prime\prime},z^\prime),x>=<y^{\prime\prime},m^*(z^\prime,x)>$ ~where $x\in X$, $y^{\prime\prime}\in Y^{**}$, $z^\prime\in Z^*$,\\
3. $m^{***}:X^{**}\times Y^{**}\rightarrow Z^{**}$,~ given by~ ~ ~$<m^{***}(x^{\prime\prime},y^{\prime\prime}),z^\prime>$ \\$=<x^{\prime\prime},m^{**}(y^{\prime\prime},z^\prime)>$ ~where ~$x^{\prime\prime}\in X^{**}$, $y^{\prime\prime}\in Y^{**}$, $z^\prime\in Z^*$.\\
The mapping $m^{***}$ is the unique extension of $m$ such that $x^{\prime\prime}\rightarrow m^{***}(x^{\prime\prime},y^{\prime\prime})$ from $X^{**}$ into $Z^{**}$ is $weak^*-weak^*$ continuous for every $y^{\prime\prime}\in Y^{**}$, but the mapping $y^{\prime\prime}\rightarrow m^{***}(x^{\prime\prime},y^{\prime\prime})$ is not in general $weak^*-weak^*$ continuous from $Y^{**}$ into $Z^{**}$ unless $x^{\prime\prime}\in X$. Hence the first topological center of $m$ may  be defined as following
$$Z_1(m)=\{x^{\prime\prime}\in X^{**}:~~y^{\prime\prime}\rightarrow m^{***}(x^{\prime\prime},y^{\prime\prime})~~is~~weak^*-weak^*~~continuous\}.$$
Let now $m^t:Y\times X\rightarrow Z$ be the transpose of $m$ defined by $m^t(y,x)=m(x,y)$ for every $x\in X$ and $y\in Y$. Then $m^t$ is a continuous bilinear map from $Y\times X$ to $Z$, and so it may be extended as above to $m^{t***}:Y^{**}\times X^{**}\rightarrow Z^{**}$.
 The mapping $m^{t***t}:X^{**}\times Y^{**}\rightarrow Z^{**}$ in general is not equal to $m^{***}$, see [1], if $m^{***}=m^{t***t}$, then $m$ is called Arens regular. The mapping $y^{\prime\prime}\rightarrow m^{t***t}(x^{\prime\prime},y^{\prime\prime})$ is $weak^*-weak^*$ continuous for every $y^{\prime\prime}\in Y^{**}$, but the mapping $x^{\prime\prime}\rightarrow m^{t***t}(x^{\prime\prime},y^{\prime\prime})$ from $X^{**}$ into $Z^{**}$ is not in general  $weak^*-weak^*$ continuous for every $y^{\prime\prime}\in Y^{**}$. So we define the second topological center of $m$ as
$$Z_2(m)=\{y^{\prime\prime}\in Y^{**}:~~x^{\prime\prime}\rightarrow m^{t***t}(x^{\prime\prime},y^{\prime\prime})~~is~~weak^*-weak^*~~continuous\}.$$
It is clear that $m$ is Arens regular if and only if $Z_1(m)=X^{**}$ or $Z_2(m)=Y^{**}$. Arens regularity of $m$ is equivalent to the following
$$\lim_i\lim_j<z^\prime,m(x_i,y_j)>=\lim_j\lim_i<z^\prime,m(x_i,y_j)>,$$
whenever both limits exist for all bounded sequences $(x_i)_i\subseteq X$ , $(y_i)_i\subseteq Y$ and $z^\prime\in Z^*$, see [18].\\
The mapping $m$ is left strongly Arens irregular if $Z_1(m)=X$ and $m$ is right strongly Arens irregular if $Z_2(m)=Y$.\\
Let now $B$ be a Banach $A-bimodule$, and let\\
$$\pi_\ell:~A\times B\rightarrow B~~~and~~~\pi_r:~B\times A\rightarrow B.$$
be the left and right module actions of $A$ on $B$, respectively. Then $B^{**}$ is a Banach $A^{**}-bimodule$ with module actions
$$\pi_\ell^{***}:~A^{**}\times B^{**}\rightarrow B^{**}~~~and~~~\pi_r^{***}:~B^{**}\times A^{**}\rightarrow B^{**}.$$
Similarly, $B^{**}$ is a Banach $A^{**}-bimodule$ with module actions\\
$$\pi_\ell^{t***t}:~A^{**}\times B^{**}\rightarrow B^{**}~~~and~~~\pi_r^{t***t}:~B^{**}\times A^{**}\rightarrow B^{**}.$$
We may therefore define the topological centers of the left and right module actions of $A$ on $B$ as follows:\\
$$Z_{B^{**}}(A^{**})=Z(\pi_\ell)=\{a^{\prime\prime}\in A^{**}:~the~map~~b^{\prime\prime}\rightarrow \pi_\ell^{***}(a^{\prime\prime}, b^{\prime\prime})~:~B^{**}\rightarrow B^{**}$$$$~is~~~weak^*-weak^*~continuous\}$$
$$Z_{B^{**}}^t(A^{**})=Z(\pi_r^t)=\{a^{\prime\prime}\in A^{**}:~the~map~~b^{\prime\prime}\rightarrow \pi_r^{t***}(a^{\prime\prime}, b^{\prime\prime})~:~B^{**}\rightarrow B^{**}$$$$~is~~~weak^*-weak^*~continuous\}$$
$$Z_{A^{**}}(B^{**})=Z(\pi_r)=\{b^{\prime\prime}\in B^{**}:~the~map~~a^{\prime\prime}\rightarrow \pi_r^{***}(b^{\prime\prime}, a^{\prime\prime})~:~A^{**}\rightarrow B^{**}$$$$~is~~~weak^*-weak^*~continuous\}$$
$$Z_{A^{**}}^t(B^{**})=Z(\pi_\ell^t)=\{b^{\prime\prime}\in B^{**}:~the~map~~a^{\prime\prime}\rightarrow \pi_\ell^{t***}(b^{\prime\prime}, a^{\prime\prime})~:~A^{**}\rightarrow B^{**}$$$$~is~~~weak^*-weak^*~continuous\}$$

\noindent We note also that if $B$ is a left(resp. right) Banach $A-module$ and $\pi_\ell:~A\times B\rightarrow B$~(resp. $\pi_r:~B\times A\rightarrow B$) is left (resp. right) module action of $A$ on $B$, then $B^*$ is a right (resp. left) Banach $A-module$. \\
We write $ab=\pi_\ell(a,b)$, $ba=\pi_r(b,a)$, $\pi_\ell(a_1a_2,b)=\pi_\ell(a_1,a_2b)$, \\ $\pi_r(b,a_1a_2)=\pi_r(ba_1,a_2)$,~
$\pi_\ell^*(a_1b^\prime, a_2)=\pi_\ell^*(b^\prime, a_2a_1)$,~
$\pi_r^*(b^\prime a, b)=\pi_r^*(b^\prime, ab)$,~ \\for all $a_1,a_2, a\in A$, $b\in B$ and  $b^\prime\in B^*$
when there is no confusion.\\
Regarding $A$ as a Banach $A-bimodule$, the operation $\pi:A\times A\rightarrow A$ extends to $\pi^{***}$ and $\pi^{t***t}$ defined on $A^{**}\times A^{**}$. These extensions are known, respectively, as the first(left) and the second (right) Arens products, and with each of them, the second dual space $A^{**}$ becomes a Banach algebra. In this situation, we shall also simplify our notations. So the first (left) Arens product of $a^{\prime\prime},b^{\prime\prime}\in A^{**}$ shall be simply indicated by $a^{\prime\prime}b^{\prime\prime}$ and defined by the three steps:
 $$<a^\prime a,b>=<a^\prime ,ab>,$$
  $$<a^{\prime\prime} a^\prime,a>=<a^{\prime\prime}, a^\prime a>,$$
  $$<a^{\prime\prime}b^{\prime\prime},a^\prime>=<a^{\prime\prime},b^{\prime\prime}a^\prime>.$$
 for every $a,b\in A$ and $a^\prime\in A^*$. Similarly, the second (right) Arens product of $a^{\prime\prime},b^{\prime\prime}\in A^{**}$ shall be  indicated by $a^{\prime\prime}ob^{\prime\prime}$ and defined by :
 $$<a oa^\prime ,b>=<a^\prime ,ba>,$$
  $$<a^\prime oa^{\prime\prime} ,a>=<a^{\prime\prime},a oa^\prime >,$$
  $$<a^{\prime\prime}ob^{\prime\prime},a^\prime>=<b^{\prime\prime},a^\prime ob^{\prime\prime}>.$$
  for all $a,b\in A$ and $a^\prime\in A^*$.\\
  The regularity of a normed algebra $A$ is defined to be the regularity of its algebra multiplication when considered as a bilinear mapping. Let $a^{\prime\prime}$ and $b^{\prime\prime}$ be elements of $A^{**}$, the second dual of $A$. By $Goldstine^,s$ Theorem [4, P.424-425], there are nets $(a_{\alpha})_{\alpha}$ and $(b_{\beta})_{\beta}$ in $A$ such that $a^{\prime\prime}=weak^*-\lim_{\alpha}a_{\alpha}$ ~and~  $b^{\prime\prime}=weak^*-\lim_{\beta}b_{\beta}$. So it is easy to see that for all $a^\prime\in A^*$,
$$\lim_{\alpha}\lim_{\beta}<a^\prime,\pi(a_{\alpha},b_{\beta})>=<a^{\prime\prime}b^{\prime\prime},a^\prime>$$ and
$$\lim_{\beta}\lim_{\alpha}<a^\prime,\pi(a_{\alpha},b_{\beta})>=<a^{\prime\prime}ob^{\prime\prime},a^\prime>,$$
where $a^{\prime\prime}b^{\prime\prime}$ and $a^{\prime\prime}ob^{\prime\prime}$ are the first and second Arens products of $A^{**}$, respectively, see [14, 18].\\
  We find the usual first and second topological center of $A^{**}$, which are
  $$Z_{A^{**}}(A^{**})=Z(\pi)=\{a^{\prime\prime}\in A^{**}: b^{\prime\prime}\rightarrow a^{\prime\prime}b^{\prime\prime}~ is~weak^*-weak^*$$$$~continuous\},$$
   $$Z^t_{A^{**}}(A^{**})=Z(\pi^t)=\{a^{\prime\prime}\in A^{**}: a^{\prime\prime}\rightarrow a^{\prime\prime}ob^{\prime\prime}~ is~weak^*-weak^*$$$$~continuous\}.$$
 An element $e^{\prime\prime}$ of $A^{**}$ is said to be a mixed unit if $e^{\prime\prime}$ is a
right unit for the first Arens multiplication and a left unit for
the second Arens multiplication. That is, $e^{\prime\prime}$ is a mixed unit if
and only if,
for each $a^{\prime\prime}\in A^{**}$, $a^{\prime\prime}e^{\prime\prime}=e^{\prime\prime}o a^{\prime\prime}=a^{\prime\prime}$. By
[4, p.146], an element $e^{\prime\prime}$ of $A^{**}$  is  mixed
      unit if and only if it is a $weak^*$ cluster point of some BAI $(e_\alpha)_{\alpha \in I}$  in
      $A$.\\
A functional $a^\prime$ in $A^*$ is said to be $wap$ (weakly almost
 periodic) on $A$ if the mapping $a\rightarrow a^\prime a$ from $A$ into
 $A^{*}$ is weakly compact. Pym in [18] showed that  this definition to the equivalent following condition\\
 For any two net $(a_{\alpha})_{\alpha}$ and $(b_{\beta})_{\beta}$
 in $\{a\in A:~\parallel a\parallel\leq 1\}$, we have\\
$$\lim_{\alpha}\lim_{\beta}<a^\prime,a_{\alpha}b_{\beta}>=\lim_{\beta}\lim_{\alpha}<a^\prime,a_{\alpha}b_{\beta}>,$$
whenever both iterated limits exist. The collection of all $wap$
functionals on $A$ is denoted by $wap(A)$. Also we have
$a^{\prime}\in wap(A)$ if and only if $<a^{\prime\prime}b^{\prime\prime},a^\prime>=<a^{\prime\prime}ob^{\prime\prime},a^\prime>$ for every $a^{\prime\prime},~b^{\prime\prime} \in
A^{**}$. \\
This paper is organized as follows:\\
\noindent {\bf a}) Let  $B$ be a left Banach  $A-bimodule$ and $\phi\in (A^{(n-r)}A^{(r)})^{(r)}$ where $0\leq r\leq n$. Then $\phi\in{ {Z}^\ell}_{B^{(n)}}((A^{(n-r)}A^{(r)})^{(r)})$ if and only if $b^{(n-1)}\phi\in B^{(n-1)}$ for all $b^{(n-1)}\in B^{(n-1)}$.\\
\noindent {\bf b}) Let  $B$ be a  Banach  $A-bimodule$. Then we have the following assertions.
\begin{enumerate}
\item ~~$b^{(n)}\in { {Z}^\ell}_{A^{(n)}}(B^{(n)})$ if and only if $b^{(n-1)}b^{(n)}\in A^{(n-1)}$ for all $b^{(n-1)}\in B^{(n-1)}$.
\item ~~ If $b^{(n-1)}b^{(n)}\in A^{(n-1)}A^*$ for each $b^{(n-1)}\in B^{(n-1)}$, then we have $b^{(n)}\in { {Z}^\ell}_{A^{(n-1)}A^*}(B^{(n)})$.
\item ~~If $\phi\in{ {Z}^\ell}_{B^{(n)}}((A^{(n-r)}A^{(r)})^{(r)})$,  then $a^{(n-1)}\phi\in { {Z}^\ell}_{B^{(n)}}(A^{(n)})$ for all ${a^{(n-1)}}\in {A^{(n-1)}}$.

\end{enumerate}

\noindent {\bf c}) Let $B$ be a Banach space such that $B^{(n)}$ is weakly compact.  Then for Banach  $A-bimodule$ $B$, we have the following assertions.

\begin{enumerate}
\item  Suppose that $(e_\alpha^{(n)})_\alpha\subseteq A^{(n)}$ is a $BLAI$ for $B^{(n)}$ such that $e_\alpha^{(n)}B^{(n+2)}\subseteq B^{(n)}$ for every $\alpha$. Then $B$ is reflexive.
\item Suppose that $(e_\alpha^{(n)})_\alpha\subseteq A^{(n)}$ is a $BRAI$ for $B^{(n)}$ and ${{Z}}^{\ell}_{e^{(n+2)}}(B^{(n+2)})=B^{(n+2)}$ where $e_\alpha^{(n)} \stackrel{w^*} {\rightarrow} e^{(n+2)}$ on $A^{(n)}$. If $B^{(n+2)}e_\alpha^{(n)}\subseteq B^{(n)}$ for every $\alpha$, then
${{Z}}^{\ell}_{A^{(n+2)}}(B^{(n+2)})=B^{(n+2)}$.
\end{enumerate}

\noindent {\bf d}) Assume that   $B$ is a  Banach  $A-bimodule$. Then we have the following assertions.
\begin{enumerate}
\item $B^{(n+1)}A^{(n)}\subseteq wap_\ell(B^{(n)})$ if and only if $A^{(n)}A^{(n+2)}\subseteq {{Z}}^{\ell}_{B^{(n+2)}}(A^{(n+2)})$.
\item If $A^{(n)}A^{(n+2)}\subseteq A^{(n)}{ {Z}^\ell}_{B^{(n+2)}}(A^{(n+2)})$, then
$$A^{(n)}A^{(n+2)}\subseteq {{Z}}^{\ell}_{B^{(n+2)}}(A^{(n+2)}).$$
\end{enumerate}

\noindent {\bf e}) Let  $B$ be a left Banach  $A-bimodule$ and $n\geq 0$ be a even. Suppose that
$b_0^{(n+1)}\in B^{(n+1)}$. Then $b_0^{(n+1)}\in wap_\ell(B^{(n)})$ if and only if the mapping $T:b^{(n+2)}\rightarrow b^{(n+2)}b_0^{(n+1)}$ form $B^{(n+2)}$ into $A^{(n+1)}$ is $weak^*-to-weak$ continuous.\\
\noindent {\bf f}) Let   $B$ be a  left Banach  $A-bimodule$. Then for $n\geq 2$, we have the following assertions.
\begin{enumerate}
\item  If $A^{(n)}=a^{(n-2)}_0A^{(n)}$ [resp.  $A^{(n)}=A^{(n)}a^{(n-2)}_0$] for some $a^{(n-2)}_0\in A^{(n-2)}$ and $a^{(n-2)}_0$ has $Rw^*w-$ property [resp. $Lw^*w-$ property] with respect to $B^{(n)}$, then $Z_{B^{(n)}}(A^{(n)})=A^{(n)}$.\\

\item  If $B^{(n)}=a^{(n-2)}_0B^{(n)}$ [resp.  $B^{(n)}=B^{(n)}a^{(n-2)}_0$] for some $a^{(n-2)}_0\in A^{(n-2)}$ and $a^{(n-2)}_0$ has $Rw^*w-$ property [resp. $Lw^*w-$ property] with respect to $B^{(n)}$, then $Z_{A^{(n)}}(B^{(n)})=B^{(n)}$.

\end{enumerate}

\begin{center}
\section{ \bf Topological centers of module actions }
\end{center}
\noindent Suppose that $A$ is a Banach algebra and $B$ is a Banach $A-bimodule$. According to [5, pp.27 and 28], $B^{**}$ is a Banach $A^{**}-bimodule$, where  $A^{**}$ is equipped with the first Arens product. So we recalled the topological centers  of module actions as in the following.

$${Z}^\ell_{A^{**}}(B^{**})=\{b^{\prime\prime}\in B^{**}:~the~map~~a^{\prime\prime}\rightarrow b^{\prime\prime} a^{\prime\prime}~:~A^{**}\rightarrow B^{**}$$$$~is~~~weak^*-weak^*~continuous\}$$
$${Z}^\ell_{B^{**}}(A^{**})=\{a^{\prime\prime}\in A^{**}:~the~map~~b^{\prime\prime}\rightarrow a^{\prime\prime} b^{\prime\prime}~:~B^{**}\rightarrow B^{**}$$$$~is~~~weak^*-weak^*~continuous\}.$$\\

\noindent Let $A^{(n)}$ and  $B^{(n)}$  be $n-th~dual$ of $B$ and $A$, respectively. By [24], $B^{(n)}$ is a Banach $A^{(n)}-bimodule$ and  we define  $B^{(n)}B^{(n-1)}$ as a subspace of $A^{(n)}$, that is, for all $b^{(n)}\in B^{(n)}$ and $b^{(n-1)}\in B^{(n-1)}$, we define
$$<b^{(n)}b^{(n-1)},a^{(n-1)}>=<b^{(n)},b^{(n-1)}a^{(n-1)}>;$$
and if $n=0$, we take $A^{(0)}=A$ and $B^{(0)}=B$.\\

\noindent {\it{\bf Theorem 2-1.}} Let  $B$ be a left Banach  $A-bimodule$ and $\phi\in (A^{(n-r)}A^{(r)})^{(r)}$ where $0\leq r\leq n$. Then $\phi\in{ {Z}^\ell}_{B^{(n)}}((A^{(n-r)}A^{(r)})^{(r)})$ if and only if $b^{(n-1)}\phi\in B^{(n-1)}$ for all $b^{(n-1)}\in B^{(n-1)}$.
\begin{proof} Let  $\phi\in{ {Z}^\ell}_{B^{(n)}}((A^{(n-r)}A^{(r)})^{(r)})$. Suppose that $(b_\alpha^{(n)})_\alpha\subseteq B^{(n)}$ such that  $b_\alpha^{(n)} \stackrel{w^*} {\rightarrow} b^{(n)}$ on $B^{(n)}$. Then, for every $b^{(n-1)}\in B^{(n-1)}$, we have
$$<b^{(n-1)}\phi,b_\alpha^{(n)}>=<b^{(n-1)},\phi b_\alpha^{(n)}>=<\phi b_\alpha^{(n)},b^{(n-1)}>
\rightarrow <\phi b^{(n)},b^{(n-1)}>$$$$=<b^{(n-1)}\phi,b^{(n)}>.$$
It follows that $b^{(n-1)}\phi\in (B^{(n+1)}, weak^*)^*=B^{(n-1)}.$\\
Conversely, let $b^{(n-1)}\phi\in B^{(n-1)}$ for every $b^{(n-1)}\in B^{(n-1)}$ and suppose that $(b_\alpha^{(n)})_\alpha\subseteq B^{(n)}$ such that  $b_\alpha^{(n)} \stackrel{w^*} {\rightarrow} b^{(n)}$ on $B^{(n)}$. Then
$$<\phi b_\alpha^{(n)},b^{(n-1)}>=<\phi, b_\alpha^{(n)}b^{(n-1)}>=< b_\alpha^{(n)}b^{(n-1)},\phi>=
< b_\alpha^{(n)},b^{(n-1)}\phi>$$$$\rightarrow < b^{(n)},b^{(n-1)}\phi>=< \phi b^{(n)},b^{(n-1)}>.$$
It follows that $<\phi b_\alpha^{(n)}\stackrel{w^*} {\rightarrow}\phi b^{(n)}$, and so $\phi\in{ {Z}^\ell}_{B^{(n)}}((A^{(n-r)}A^{(r)})^{(r)})$.
\end{proof}

\noindent Let  $B$ be a left Banach  $A-bimodule$ and $\phi\in (A^{(n-r)}A^{(r)})^{(n-r)}$ where $0\leq r\leq n$. Then $\phi\in{ {Z}^\ell}_{B^{(n)}}((A^{(n-r)}A^{(r)})^{(n-r)})$ if and only if $b^{(n-1)}\phi\in B^{(n-1)}$ for every $b^{(n-1)}\in B^{(n-1)}$. The proof of the proceeding  assertion is the similar to Theorem 2-1, and
if we take $B=A$,  $n=1$ and $r=0$, we obtain Lemma 3.1 (b) from [14].\\

\noindent {\it{\bf Theorem 2-2.}} Let  $B$ be a  Banach  $A-bimodule$. Then we have the following assertions.
\begin{enumerate}
\item ~~$b^{(n)}\in { {Z}^\ell}_{A^{(n)}}(B^{(n)})$ if and only if $b^{(n-1)}b^{(n)}\in A^{(n-1)}$ for all $b^{(n-1)}\in B^{(n-1)}$.
\item ~~ If $b^{(n-1)}b^{(n)}\in A^{(n-1)}A^*$ for each $b^{(n-1)}\in B^{(n-1)}$, then we have $b^{(n)}\in { {Z}^\ell}_{A^{(n-1)}A^*}(B^{(n)})$.
\item ~~If $\phi\in{ {Z}^\ell}_{B^{(n)}}((A^{(n-r)}A^{(r)})^{(r)})$,  then $a^{(n-1)}\phi\in { {Z}^\ell}_{B^{(n)}}(A^{(n)})$ for all ${a^{(n-1)}}\in {A^{(n-1)}}$.

\end{enumerate}
\begin{proof}
\begin{enumerate}
\item  Let ~$b^{(n)}\in { {Z}^\ell}_{A^{(n)}}(B^{(n)})$. We show that $b^{(n-1)}b^{(n)}\in A^{(n-1)}$ where $b^{(n-1)}\in B^{(n-1)}$. Suppose that
$(a_\alpha^{(n)})_\alpha\subseteq A^{(n)}$ and  $a_\alpha^{(n)} \stackrel{w^*} {\rightarrow} a^{(n)}$ on $A^{(n)}$.     Then we have
$$<b^{(n-1)}b^{(n)},a_\alpha^{(n)}>=<b^{(n-1)},b^{(n)}a_\alpha^{(n)}>=<b^{(n)}a_\alpha^{(n)},b^{(n-1)}>$$
$$\rightarrow<b^{(n)}a^{(n)},b^{(n-1)}>=
 <b^{(n-1)}b^{(n)},a^{(n)}>.$$
 Consequently  $b^{(n-1)}b^{(n)}\in (A^{(n+1)},weak^*)^*=A^{(n-1)}$.  It follows that \\ $b^{(n-1)}b^{(n)}\in A^{(n-1)}$.\\
 Conversely, let $b^{(n-1)}b^{(n)}\in A^{(n-1)}$ for each $b^{(n-1)}\in B^{(n-1)}$. Suppose that
$(a_\alpha^{(n)})_\alpha\subseteq A^{(n)}$ and  $a_\alpha^{(n)} \stackrel{w^*} {\rightarrow} a^{(n)}$ on $A^{(n)}$.   Then we have
 $$<b^{(n)}a_\alpha^{(n)},b^{(n-1)}>=<b^{(n)},a_\alpha^{(n)}b^{(n-1)}>= <a_\alpha^{(n)}b^{(n-1)},b^{(n)}>$$$$=<a_\alpha^{(n)},b^{(n-1)}b^{(n)}>\rightarrow <a^{(n)},b^{(n-1)}b^{(n)}>=
 <b^{(n)}a^{(n)},b^{(n-1)}>.$$
It follows that  $b^{(n)}a_\alpha^{(n)} \stackrel{w^*}\rightarrow b^{(n)}a^{(n)}$, and so  $b^{(n)}\in { {Z}^\ell}_{A^{(n)}}(B^{(n)})$.
\item Proof is similar to (1).
\item Let $\phi\in{ {Z}^\ell}_{B^{(n)}}((A^{(n-r)}A^{(r)})^{(r)})$ and $a^{(n-1)}\in A^{(n-1)}$. Assume that
$(b_\alpha^{(n)})_\alpha\subseteq B^{(n)}$ and  $b_\alpha^{(n)} \stackrel{w^*} {\rightarrow} b^{(n)}$ on $B^{(n)}$. Then for all $b^{(n-1)}\in B^{(n-1)}$, we have
$$<(a^{(n-1)}\phi) b_\alpha^{(n)},b^{(n-1)}>=<\phi b_\alpha^{(n)},b^{(n-1)}a^{(n-1)}>\rightarrow
<\phi b^{(n)},b^{(n-1)}a^{(n-1)}>$$$$=<(a^{(n-1)}\phi) b^{(n)},b^{(n-1)}>.$$
It follows that $(a^{(n-1)}\phi) b_\alpha^{(n)}\stackrel{w^*} {\rightarrow}(a^{(n-1)}\phi) b^{(n)}$, and so
$a^{(n-1)}\phi\in { {Z}^\ell}_{B^{(n)}}(A^{(n)})$.
\end{enumerate}
\end{proof}
In the proceeding theorem, part (1), if we take $B=A$ and $n=2$, we conclude Lemma 3.1 (a) from [14].  We change part (3) of Theorem 2-2 with the following statement, that is, if $\phi\in{ {Z}^\ell}_{B^{(n)}}((A^{(n-r)}A^{(r)})^{(n-r)})$,  then $a^{(n-1)}\phi\in { {Z}^\ell}_{B^{(n)}}(A^{(n)})$ for all ${a^{(n-1)}}\in {A^{(n-1)}}$. Proof of this claim is similar to proceeding theorem and if we take $B=A$, $n=1$ and $r=0$,  we obtain Lemma 3.1 (c) from [14].\\

\noindent {\it{\bf Definition 2-3.}} Let $B$ be a Banach  $A-bimodule$ and suppose that  $a^{\prime\prime}\in A^{**}$. We say that $a^{\prime\prime}\rightarrow b^{\prime\prime}a^{\prime\prime}$~is~$weak^*-to-weak^*$~point~~continuous, if for every net $(a_\alpha^{\prime\prime})_\alpha\subseteq A^{**}$ such that
$a^{\prime\prime}_\alpha\stackrel{w^*} {\rightarrow}a^{\prime\prime}$, it follows that $a^{\prime\prime}_\alpha b^{\prime\prime}   \stackrel{w^*} {\rightarrow}a^{\prime\prime} b^{\prime\prime}$.\\
 Suppose that  $B$ is a   Banach $A-bimodule$. Assume that $a^{\prime\prime}\in A^{**}$. Then we define  the  locally topological center of $a^{\prime\prime}$ on $B^{**}$ as follows
$${{Z}}^\ell_{a^{\prime\prime}}(B^{**})=\{ b^{\prime\prime}\in B^{**}:~a^{\prime\prime}\rightarrow b^{\prime\prime}a^{\prime\prime}~is~weak^*-to-weak^*~point~~continuous\},$$

The definition of  ${ {Z}}^\ell_{b^{\prime\prime}}(A^{**})$ where
$b^{\prime\prime}\in B^{**}$ are similar.\\
It is clear that $$\bigcap_{a^{\prime\prime}\in A^{**}}{{Z}}^{\ell}_{a^{\prime\prime}}(B^{**})={{Z}}^{\ell}_{A^{**}}(B^{**}),$$
$$\bigcap_{b^{\prime\prime}\in B^{**}}{{Z}}^{\ell}_{b^{\prime\prime}}(A^{**})={{Z}}^{\ell}_{B^{**}}(A^{**}).$$\\

\noindent {\it{\bf Theorem 2-4.}} Let $B$ be a Banach space such that $B^{(n)}$ is weakly compact.  Then for Banach  $A-bimodule$ $B$, we have the following assertions.

\begin{enumerate}
\item  Suppose that $(e_\alpha^{(n)})_\alpha\subseteq A^{(n)}$ is a $BLAI$ for $B^{(n)}$ such that $e_\alpha^{(n)}B^{(n+2)}\subseteq B^{(n)}$ for every $\alpha$. Then $B$ is reflexive.
\item Suppose that $(e_\alpha^{(n)})_\alpha\subseteq A^{(n)}$ is a $BRAI$ for $B^{(n)}$ and ${{Z}}^{\ell}_{e^{(n+2)}}(B^{(n+2)})=B^{(n+2)}$ where $e_\alpha^{(n)} \stackrel{w^*} {\rightarrow} e^{(n+2)}$ on $A^{(n)}$. If $B^{(n+2)}e_\alpha^{(n)}\subseteq B^{(n)}$ for every $\alpha$, then
${{Z}}^{\ell}_{A^{(n+2)}}(B^{(n+2)})=B^{(n+2)}$.
\end{enumerate}

\begin{proof}
\begin{enumerate}
\item  Let $b^{n+2}\in B^{n+2}$. Since $(e_\alpha^{(n)})_\alpha$ is a $BLAI$ for $B^{(n)}$, there is left unit as  $e^{(n+2)}\in A^{n+2}$ for $B^{n+2}$, see [10]. Then we have  $e_\alpha^{(n)} b^{(n+2)} \stackrel{w^*} {\rightarrow} b^{(n+2)}$ on $B^{(n+2)}$. Since $e_\alpha^{(n)}b^{(n+2)}\in B^{(n)}$, we have  $e_\alpha^{(n)} b^{(n+2)} \stackrel{w} {\rightarrow} b^{(n+2)}$ on $B^{(n)}$. We conclude that $b^{n+2}\in B^{n}$ of course $B^{n}$ is weakly compact.
\item Suppose that $b^{(n+2)}\in{{Z}}^{\ell}_{A^{(n+2)}}(B^{(n+2)})$ and  $e_\alpha^{(n)} \stackrel{w^*} {\rightarrow} e^{(n+2)}$ on $A^{(n)}$ such that $e^{(n+2)}$ is right unit for $B^{(n+2)}$, see [10]. Then we have
$ b^{(n+2)}e_\alpha^{(n)} \stackrel{w^*} {\rightarrow} b^{(n+2)}$ on $B^{(n+2)}$. Since $B^{(n+2)}e_\alpha^{(n)}\subseteq B^{(n)}$ for every $\alpha$, $ b^{(n+2)}e_\alpha^{(n)} \stackrel{w} {\rightarrow} b^{(n+2)}$ on $B^{(n)}$ and since  $B^{(n)}$ is weakly compact, $b^{(n+2)}\in B^{(n)}$. It follows that
${{Z}}^{\ell}_{A^{(n+2)}}(B^{(n+2)})=B^{(n+2)}$.
\end{enumerate}\end{proof}

\noindent {\it{\bf Definition 2-5.}} Let  $B$ be a  Banach  $A-bimodule$. Then $b^{n+1}\in B^{n+1}$ is said to be weakly left almost periodic functional if the set
$$\{b^{(n+1)}a^{(n)}:~a^{(n)}\in A^{(n)}, ~ \parallel a^{(n)}\parallel\leq 1\},$$
is relatively weakly compact, and $b^{n+1}\in B^{n+1}$ is said to be weakly right almost periodic functional if the set
$$\{a^{(n)}b^{(n+1)}:~a^{(n)}\in A^{(n)}, ~ \parallel a^{(n)}\parallel\leq 1\},$$
is relatively weakly compact. We denote by $wap_\ell(B^{(n)})$ [resp. $wap_\ell(B^{(n)})$] the closed subspace of $B^{(n+1)}$ consisting of all the weakly left [resp. right] almost periodic functionals in $B^{(n+1)}$. By [6, 14, 18], the definition of $wap_\ell(B^{(n)})$ and  $wap_r(B^{(n)})$, respectively, are  equivalent to the following
$$wap_\ell(B^{(n)})=\{b^{(n+1)}\in B^{(n+1)}:~<b^{(n+2)}a_\alpha^{(n+2)},b^{(n+1)}>\rightarrow <b^{(n+2)}a^{(n+2)},b^{(n+1)}> $$$$~~where~~ a_{\alpha}^{(n+2)}  \stackrel{w^*} {\rightarrow} a^{(n+2)} \}.$$
and
$$wap_r(B^{(n)})=\{b^{(n+1)}\in B^{(n+1)}:~<a^{(n+2)}b_\alpha^{(n+2)},b^{(n+1)}>\rightarrow <a^{(n+2)}b^{(n+2)},b^{(n+1)}> $$$$~~where~~ b_{\alpha}^{(n+2)}  \stackrel{w^*} {\rightarrow} b^{(n+2)} \}.$$
If we take $A=B$ and $n=0$, then $wap_\ell(A)=wap_r(A)=wap(A)$.\\

\noindent {\it{\bf Theorem 2-6.}} Assume that   $B$ is a  Banach  $A-bimodule$. Then we have the following assertions.
\begin{enumerate}
\item $B^{(n+1)}A^{(n)}\subseteq wap_\ell(B^{(n)})$ if and only if $A^{(n)}A^{(n+2)}\subseteq {{Z}}^{\ell}_{B^{(n+2)}}(A^{(n+2)})$.
\item If $A^{(n)}A^{(n+2)}\subseteq A^{(n)}{ {Z}^\ell}_{B^{(n+2)}}(A^{(n+2)})$, then
$$A^{(n)}A^{(n+2)}\subseteq {{Z}}^{\ell}_{B^{(n+2)}}(A^{(n+2)}).$$
\end{enumerate}
\begin{proof}
\begin{enumerate}
\item Suppose that $B^{(n+1)}A^{(n)}\subseteq wap_\ell(B^{(n)})$. Let  $a^{(n)}\in A^{(n)}$, $a^{(n+2)}\in A^{(n+2)}$ and let $(b_{\alpha}^{(n+2)})_{\alpha}\subseteq B^{(n+2)}$ such that   $ b_{\alpha}^{(n+2)}  \stackrel{w^*} {\rightarrow} b^{(n+2)}$. Then for every $b^{(n+1)}\in B^{(n+1)}$, we have
$$<(a^{(n)}a^{(n+2)})b_{\alpha}^{(n+2)},b^{(n+1)}>=<a^{(n+2)}b_{\alpha}^{(n+2)},b^{(n+1)}a^{(n)}>$$$$\rightarrow
<a^{(n+2)}b^{(n+2)},b^{(n+1)}a^{(n)}>=<(a^{(n)}a^{(n+2)})b^{(n+2)},b^{(n+1)}>.$$
It follows that $a^{(n)}a^{(n+2)}\in {{Z}}^{\ell}_{B^{(n+2)}}(A^{(n+2)})$.\\
Conversely, let $a^{(n)}a^{(n+2)}\in {{Z}}^{\ell}_{B^{(n+2)}}(A^{(n+2)})$ for every $a^{(n)}\in A^{(n)}$, $a^{(n+2)}\in A^{(n+2)}$ and suppose that $(b_{\alpha}^{(n+2)})_{\alpha}\subseteq B^{(n+2)}$ such that   $ b_{\alpha}^{(n+2)}  \stackrel{w^*} {\rightarrow} b^{(n+2)}$. Then for every $b^{(n+1)}\in B^{(n+1)}$, we have
$$<a^{(n+2)}b_{\alpha}^{(n+2)},b^{(n+1)}a^{(n)}>=<(a^{(n)}a^{(n+2)})b_{\alpha}^{(n+2)},b^{(n+1)}>$$$$\rightarrow
<(a^{(n)}a^{(n+2)})b^{(n+2)},b^{(n+1)}>=<a^{(n+2)}b_{\alpha}^{(n+2)},b^{(n+1)}a^{(n)}>.$$
It follows that $B^{(n+1)}A^{(n)}\subseteq wap_\ell(B^{(n)})$.
\item Since $A^{(n)}A^{(n+2)}\subseteq A^{(n)}{ {Z}^\ell}_{B^{(n)}}((A^{(n+2)})$, for every $a^{(n)}\in A^{(n)}$ and  $a^{(n+2)}\in A^{(n+2)}$, we have $a^{(n)}a^{(n+2)}\in A^{(n)}{ {Z}^\ell}_{B^{(n+2)}}(A^{(n+2)})$.

     Then there are $x^{(n)}\in A^{(n)}$ and $\phi\in {Z}^\ell_{B^{(n+2)}}(A^{(n+2)})$ such that
   $a^{(n)}a^{(n+2)}=x^{(n)}\phi$. Suppose that $(b_{\alpha}^{(n+2)})_{\alpha}\subseteq B^{(n+2)}$ such that

     $b_{\alpha}^{(n+2)}  \stackrel{w^*} {\rightarrow} b^{(n+2)}$. Then for every $b^{(n+1)}\in B^{(n+1)}$, we have
$$<(a^{(n)}a^{(n+2)})b_{\alpha}^{(n+2)},b^{(n+1)}>=<(x^{(n)}\phi)b_{\alpha}^{(n+2)},b^{(n+1)}>$$$$=
<\phi b_{\alpha}^{(n+2)},b^{(n+1)}x^{(n)}>\rightarrow <\phi b^{(n+2)},b^{(n+1)}x^{(n)}>$$$$=
<(a^{(n)}a^{(n+2)})b^{(n+2)},b^{(n+1)}>.$$
\end{enumerate}\end{proof}

If we take $B=A$ and $n=0$, we conclude Theorem 3.6 (a) from [14].\\

\noindent {\it{\bf Theorem 2-7.}} Assume that   $B$ is a  Banach  $A-bimodule$. If $A^{(n)}$ is a left ideal in
$A^{(n+2)}$, then $B^{(n+1)}A^{(n)}\subseteq wap_\ell(B^{(n)})$.\\\\
\noindent {\it{\bf Theorem 2-8.}} Let  $B$ be a left Banach  $A-bimodule$ and $n\geq 0$ be a even. Suppose that
$b_0^{(n+1)}\in B^{(n+1)}$. Then $b_0^{(n+1)}\in wap_\ell(B^{(n)})$ if and only if the mapping $T:b^{(n+2)}\rightarrow b^{(n+2)}b_0^{(n+1)}$ form $B^{(n+2)}$ into $A^{(n+1)}$ is $weak^*-to-weak$ continuous.
\begin{proof}
Let $b_0^{(n+1)}\in B^{(n+1)}$ and suppose that  $ b^{(n+2)}_\alpha\stackrel{w^*} {\rightarrow} b^{(n+2)}$ on $B^{(n+2)}$. Then for every $a^{(n+2)}\in A^{(n+2)}$, we have
$$<a^{(n+2)},b_\alpha^{(n+2)}b_0^{(n+1)}>=<a^{(n+2)}b_\alpha^{(n+2)},b_0^{(n+1)}>\rightarrow
<a^{(n+2)}b^{(n+2)},b_0^{(n+1)}>$$$$=<a^{(n+2)},b^{(n+2)}b_0^{(n+1)}>.$$
It follows that $b_\alpha^{(n+2)}b_0^{(n+1)}\stackrel{w} {\rightarrow}b^{(n+2)}b_0^{(n+1)}$ on $A^{(n+1)}$.\\
The proof of the converse is similar of proceeding proof.
\end{proof}

\noindent {\it{\bf Corollary 2-9.}} Assume that   $B$ is a  Banach  $A-bimodule$. Then ${{Z}}^{\ell}_{A^{(n+2)}}(B^{(n+2)})=B^{(n+2)}$ if and only if the mapping $T:b^{(n+2)}\rightarrow b^{(n+2)}b_0^{(n+1)}$ form $B^{(n+2)}$ into $A^{(n+1)}$ is $weak^*-to-weak$ continuous for every $b_0^{(n+1)}\in B^{(n+1)}$.\\

\noindent {\it{\bf Corollary 2-10.}} Let $A$ be a Banach algebra.
 Assume that $a^\prime\in A^*$ and $T_{a^\prime}$ is the linear operator from
$A$ into $A^*$ defined by $T_{a^\prime} a=a^\prime a$. Then, $a^\prime\in wap(A)$ if and
only if the adjoint of $T_{a^\prime}$ is
$weak^*-to-weak$ continuous. So $A$ is Arens regular if and only if the adjoint of the mapping  $T_{a^\prime} a=a^\prime a$ is $weak^*-to-weak$ continuous for every $a^\prime\in A^*$. \\

\noindent {\it{\bf Definition 2-11.}} Let   $B$ be a  left Banach  $A-bimodule$. We say that $a^{(n)}\in A^{(n)}$ has $Left-weak^*-weak$ property ($=Lw^*w-$ property) with respect to $B^{(n)}$, if for every $(b^{(n+1)}_{\alpha})_{\alpha}\subseteq B^{(n+1)}$, $a^{(n)}b_\alpha^{(n+1)}\stackrel{w^*} {\rightarrow}0$ implies $a^{(n)}b_\alpha^{(n+1)}\stackrel{w} {\rightarrow}0$. If every $a^{(n)}\in A$ has $Lw^*w-$ property with respect to $B^{(n)}$, then we say that $A^{(n)}$ has $Lw^*w-$ property with respect to $B^{(n)}$. The definition of the $Right-weak^*-weak$ property ($=Rw^*w-$ property) is the same.\\
We say that $a^{(n)}\in A^{(n)}$ has $weak^*-weak$  property ($=w^*w-$ property) with respect to $B^{(n)}$ if it has $Lw^*w-$ property and $Rw^*w-$ property with respect to $B^{(n)}$.\\
If $a^{(n)}\in A^{(n)}$ has $Lw^*w-$ property with respect to itself, then we say that $a^{(n)}\in A^{(n)}$ has $Lw^*w-$ property.\\

\noindent {\it{\bf Example 2-12.}}
\begin{enumerate}
\item If $B$ is Banach $A$-bimodule and reflexive, then $A$ has $w^*w-$property with respect to $B$.
\item $L^1(G)$, $M(G)$ and $A(G)$ have $w^*w-$property when $G$ is finite.
\item Let $G$ be locally compact group.  $L^1(G)$ [resp. $M(G)$] has $w^*w-$property [resp. $Lw^*w-$ property ]with respect to $L^p(G)$ whenever $p>1$.
\item Suppose that $B$ is a left Banach $A-module$ and $e$ is left unit element of $A$ such that $eb=b$ for all $b\in B$. If $e$ has $Lw^*w-$ property, then $B$ is reflexive.
\item If $S$ is a compact semigroup, then $C^+(S)=\{f\in C(S):~f>0\}$ has $w^*w-$property.\\
\end{enumerate}

\noindent {\it{\bf Theorem 2-13.}} Let   $B$ be a  left Banach  $A-bimodule$. Then for $n\geq 2$, we have the following assertions.
\begin{enumerate}
\item  If $A^{(n)}=a^{(n-2)}_0A^{(n)}$ [resp.  $A^{(n)}=A^{(n)}a^{(n-2)}_0$] for some $a^{(n-2)}_0\in A^{(n-2)}$ and $a^{(n-2)}_0$ has $Rw^*w-$ property [resp. $Lw^*w-$ property] with respect to $B^{(n)}$, then $Z_{B^{(n)}}(A^{(n)})=A^{(n)}$.\\

\item  If $B^{(n)}=a^{(n-2)}_0B^{(n)}$ [resp.  $B^{(n)}=B^{(n)}a^{(n-2)}_0$] for some $a^{(n-2)}_0\in A^{(n-2)}$ and $a^{(n-2)}_0$ has $Rw^*w-$ property [resp. $Lw^*w-$ property] with respect to $B^{(n)}$, then $Z_{A^{(n)}}(B^{(n)})=B^{(n)}$.

\end{enumerate}

\begin{proof}

\begin{enumerate}
\item  Suppose that $A^{(n)}=a^{(n-2)}_0A^{(n)}$  for some $a^{(n-2)}_0\in A$ and $a^{(n-2)}_0$ has $Rw^*w-$ property. Let   $(b^{(n)}_{\alpha})_{\alpha}\subseteq B^{(n)}$ such that  $b^{(n)}_{\alpha} \stackrel{w^*} {\rightarrow}b^{(n)}$. Then for every $a^{(n-2)}\in A^{(n-2)}$ and $b^{(n-1)}\in B^{(n-1)}$, we have
$$<b_\alpha^{(n)}b^{(n-1)},a^{(n-2)}>=<b_\alpha^{(n)},b^{(n-1)}a^{(n-2)}>\rightarrow
<b^{(n)},b^{(n-1)}a^{(n-2)}>$$$$=<b^{(n)}b^{(n-1)},a^{(n-2)}>.$$
It follows that $b_\alpha^{(n)}b^{(n-1)}\stackrel{w^*} {\rightarrow} b^{(n)}b^{(n-1)}$. Also it is clear that
$(b_\alpha^{(n)}b^{(n-1)})a^{(n-2)}_0\stackrel{w^*} {\rightarrow} (b^{(n)}b^{(n-1)})a^{(n-2)}_0$.
Since $a^{(n-2)}_0$ has $Rw^*w-$ property, $(b_\alpha^{(n)}b^{(n-1)})a^{(n-2)}_0\stackrel{w} {\rightarrow} (b^{(n)}b^{(n-1)})a^{(n-2)}_0$. Now, let $a^{(n)}\in A^{(n)}$. Since $A^{(n)}=a^{(n-2)}_0A^{(n)}$, there is  $x^{(n)}\in A^{(n)}$ such that $a^{(n)}=a_0^{(n-2)}x^{(n)}$. Thus we have
$$<a^{(n)}b_\alpha^{(n)},b^{(n-1)}>=<a^{(n)},b_\alpha^{(n)}b^{(n-1)}>=<a_0^{(n-2)}x^{(n)},b_\alpha^{(n)}b^{(n-1)}>
$$$$=<x^{(n)},(b_\alpha^{(n)}b^{(n-1)})a_0^{(n-2)}>\rightarrow <x^{(n)},(b^{(n)}b^{(n-1)})a_0^{(n-2)}>$$$$=<a^{(n)}b,b^{(n-1)}>.$$
It follows that  $a^{(n)}\in Z_{A^{(n)}}(B^{(n)})$.\\
Proof of the next part is similar to proceeding proof.
\item Let $B^{(n)}=a^{(n-2)}_0B^{(n)}$  for some $a^{(n-2)}_0\in A$ and $a^{(n-2)}_0$ has $Rw^*w-$ property  with respect to $B^{(n)}$. Assume that
$(a_{\alpha}^{(n)})_{\alpha}\subseteq A^{(n)}$ such that  $a^{(n)}_{\alpha} \stackrel{w^*} {\rightarrow}a^{(n)}$. Then for every $b^{(n-1)}\in B^{(n-1)}$, we have
$$<a_\alpha^{(n)}b^{(n-1)},b^{(n-2)}>=<a_\alpha^{(n)},b^{(n-1)}b^{(n-2)}>\rightarrow <a^{(n)},b^{(n-1)}b^{(n-2)}>$$$$=
<a^{(n)}b^{(n-1)},b^{(n-2)}>$$
 We conclude that $a_\alpha^{(n)}b^{(n-1)}\stackrel{w^*} {\rightarrow} a^{(n)}b^{(n-1)}$. It is clear that
$(a_\alpha^{(n)}b^{(n-1)})a^{(n-2)}_0\stackrel{w^*} {\rightarrow} (a^{(n)}b^{(n-1)})a^{(n-2)}_0$.
Since $a^{(n-2)}_0$ has $Rw^*w-$ property, $(a_\alpha^{(n)}b^{(n-1)})a^{(n-2)}_0\stackrel{w} {\rightarrow} (a^{(n)}b^{(n-1)})a^{(n-2)}_0$.\\
Suppose that $ b^{(n)}\in B^{(n)}$. Since $B^{(n)}=a^{(n-2)}_0B^{(n)}$, there is $ y^{(n)}\in B^{(n)}$ such that
$b^{(n)}=a^{(n-2)}_0y^{(n)}$. Consequently, we have
$$<b^{(n)}a_\alpha^{(n)},b^{(n-1)}>=<b^{(n)},a_\alpha^{(n)}b^{(n-1)}>=<a^{(n-2)}_0y^{(n)},a_\alpha^{(n)}b^{(n-1)}>$$$$=
<y^{(n)},(a_\alpha^{(n)}b^{(n-1)})a^{(n-2)}_0>\rightarrow <y^{(n)},(a^{(n)}b^{(n-1)})a^{(n-2)}_0>$$$$=
<a^{(n-2)}_0y^{(n)},(a^{(n)}b^{(n-1)})>=<b^{(n)}a^{(n)},b^{(n-1)}>.$$
Thus $b^{(n)}a_\alpha^{(n)}\stackrel{w} {\rightarrow} b^{(n)}a^{(n)}$. It follows that $b^{(n)}\in Z_{A^{(n)}}(B^{(n)})$.\\
The next part similar to the proceeding proof.

\end{enumerate}\end{proof}

\noindent {\it{\bf Example 2-14.}} Let $G$ be a locally compact group. Since $M(G)$ is a Banach $L^1(G)$-bimodule and the unit element of  $M(G)^{(n)}$ has not $Lw^*w-$ property or $Rw^*w-$ property, by Theorem 2-13, $Z_{L^1(G)^{(n)}}{(M(G)^{(n)})}\neq {M(G)^{(n)}}$.\\
ii) If $G$ is finite, then by Theorem 2-13, we have $Z_{M(G)^{(n)}}{(L^1(G)^{(n)})}= {L^1(G)^{(n)}}$ and
$Z_{L^1(G)^{(n)}}{(M(G)^{(n)})}= {M(G)^{(n)}}$.\\\\

\bibliographystyle{amsplain}

\it{Department of Mathematics, Amirkabir University of Technology, Tehran, Iran\\
{\it Email address:} haghnejad@aut.ac.ir\\\\
Department of Mathematics, Amirkabir University of Technology, Tehran, Iran\\
{\it Email address:} riazi@aut.ac.ir}
\end{document}